# Soft $T_{(0,\alpha)}$ Spaces

Mbekezeli Sibahle Nxumalo

**Abstract:** We present a construction of the soft $T_0$ reflection of a soft topological space and characterize some separation axioms developed through the soft $T_0$ reflection.

## 1. Introduction

In 1999, Molodstov [10] initiated a novel concept of soft set theory which is a new approach for modelling uncertainties that are in engineering, economics, environments, medical sciences, etc. He applied the soft set theory into several directions such as smoothness of functions, game theory, Riemann integration, Perron integration, and so on.

In 2011, Shabir and Naz [12] initiated the study of soft topological spaces which are defined over an initial universe with a fixed set of parameters.

Herrlich and Strecker gave a construction of the $T_0$ Reflection on **Top** (category of topological spaces whose morphisms are continuous functions) [4]; in [5], Kunzi and Richmond constructed the $T_0$-ordered reflection on **PREORDTOP** (category of preorder topological spaces whose morphisms are continuous increasing maps); in [9], Mirhosseinkhani introduced a construction of $T_0$ Reflection on **GenTOP** (category of generalized topological spaces whose morphisms are g-continuous maps) and in [7], Abdelwaheb Mhemdi et al. considered a construction of the $T_0$ Reflection on **PreTOP** (category of pretopological spaces whose arrows are p-continuous maps). In this paper, we consider a construction of the soft $T_0$ Reflection in the category $\mathbf{STop}_{0U}$ (whose objects are soft $T_{0U}$ topological spaces and morphisms are soft continuous maps). $\mathbf{STop}_{0U}$ is a full subcategory of **STop** (category of soft topological spaces whose morphisms are soft continuous maps). We further introduce and characterize some separation axioms developed through the soft $T_0$ reflection.

## 2. Preliminaries

### 2.1  Soft Sets

**Definition 2.1.1** *[10] Let $X$ be an initial universe and $A$ be a set of parameters. Let $\mathcal{P}(X)$ denote the power set of $X$. A pair $(F, A)$, where $F$ is a map from $A$ to $\mathcal{P}(X)$, is called a soft set over $X$.*

Throughout the text, we denote the family of all soft sets $(F, A)$ over $X$ by $SS(X)_A$.

**Definition 2.1.2.** *[3] The soft set $(F, A) \in SS(X)_A$, where $F(a) = \emptyset$, for every $a \in A$ is called the A-null soft set of $SS(X)_A$ and denoted by $0_A$. The soft set $(F, A) \in SS(X)_A$, where $F(a) = X$, for every $a \in A$, is called the A-absolute soft set of $SS(X)_A$ and denoted by $1_A$. If the same set of parameters A also applies to a set Y, and X and Y appear in the same sentence, then we shall denote by $1_A^Y$, the soft set $(F, A)$ on Y such that $F(a) = Y$, for all $a \in A$.*

**Definition 2.1.3.** *[3] The soft complement of $(F, A)$ is the soft set $(H, A) \in SS(X)_A$, where the map $H : A \to \mathcal{P}(X)$ is defined as $H(a) = X \setminus F(a)$, for every $a \in A$. Symbolically, we write $(H, A) = (F, A)^c$ or $(H, A) = 1_A \setminus (F, A)$. For $(G, A) \in SS(X)_A$, $(G, A) \setminus (H, A)$ is the soft set $(F, A)$, where the map $F : A \to \mathcal{P}(X)$ is defined as $F(a) = G(a) \setminus H(a)$, for all $a \in A$.*

**Definition 2.1.4.** *[3] The soft union of soft sets $(F, A)$ and $(G, A)$ is the soft set $(H, A) \in SS(X)_A$, where the map $H : A \to \mathcal{P}(X)$ is defined as follows: $H(a) = F(a) \cup G(a)$, for every $a \in A$. Symbolically, we write $(H, A) = (F, A) \sqcup (G, A)$. The soft intersection of soft sets $(F, A)$ and $(G, A)$ is the soft set $(H, A) \in SS(X)_A$ where the map $H : A \to \mathcal{P}(X)$ is defined as follows: $H(a) = F(a) \cap G(a)$, for every $a \in A$. Symbolically, we write $(H, A) = (F, A) \sqcap (G, A)$.*

**Definition 2.1.5.** *[3] Let $(F, A), (G, A) \in SS(X)_A$. We say that $(F, A)$ is a soft subset of $(G, A)$ if $F(a) \subseteq G(a)$, for every $a \in A$. Symbolically, we write $(F, A) \sqsubseteq (G, A)$. Two soft sets $(F, A)$ and $(G, A)$ are soft equal if $(F, A) \sqsubseteq (G, A)$ and $(G, A) \sqsubseteq (F, A)$. Symbolically, we write $(F, A) = (G, A)$. Any two soft sets $(F, A)$ and $(G, A)$ are soft disjoint if $F(a) \cap G(a) = \emptyset$, for all $a \in A$. Symbolically, we write $(F, A) \sqcap (G, A) = 0_A$.*

**Definition 2.1.6.** *[12] Let $(F, A) \in SS(X)_A$ and $x \in X$. We say that $x \in (F, A)$, read as x belongs to the soft set $(F, A)$, whenever $x \in F(a)$, for all $a \in A$. For any $x \in X$, $x \notin (F, A)$, if $x \notin F(a)$, for some $a \in A$. Let $x \in X$. Then $(F_x, A)$ denotes the soft set for which $F_x(a) = \{x\}$, for all $x \in X$.*

**Definition 2.1.7.** *[12] Let $(F, A) \in SS(X)_A$ and Y be a nonempty subset of X. Then the sub soft set of $(F, A)$ over Y, denoted by $(F, A)^Y$, is a soft set $(G, A)$ defined as $G(a) = Y \cap F(a)$, for each $a \in A$.*

**Definition 2.1.8.** *[3] Let X and Y be two initial universe sets, A and B two sets of parameters, and $f : X \to Y$ and $e : A \to B$ be functions. The mapping $(f, e) : SS(X)_A \to SS(Y)_B$ is called a soft mapping from X to Y, for which:*

a. *If $(F, A) \in SS(X)_A$, then the image of $(F, A)$ under $(f, e)$, denoted by $(f, e)((F, A))$, is the soft set $(G, B) \in SS(Y)_B$ such that*

$$G(b) = \begin{cases} \cup \{f(F(a)) : a \in e^{-1}(\{b\})\}, & \text{if } e^{-1}(\{b\}) \neq \emptyset \\ \emptyset, & \text{if } e^{-1}(\{b\}) = \emptyset \end{cases}$$

      *for all* $b \in B$.
- b. *If* $(G, B) \in SS(Y)_B$, *then the pre-image of* $(G, B)$ *under* $(f, e)$, *denoted by* $(f, e)^{-1}((G, B))$, *is the soft set* $(F, A) \in SS(X)_A$ *such that* $F(a) = f^{-1}(G(e(a)))$, *for every* $a \in A$.

**Definition 2.1.9.** *[1] A soft mapping* $(f, e): SS(X)_A \to SS(Y)_B$ *is called injective, if* $f$ *and* $e$ *are injective. It is called surjective, if* $f$ *and* $e$ *are surjective.*

**Definition 2.1.10.** *Two soft mappings* $(f, e)$ *and* $(g, e')$ *are said to be equal if* $f = g$ *and* $e = e'$.

**Definition 2.1.11.** *[1] Let* $(f, e): SS(X)_A \to SS(Y)_B$ *and* $(g, e'): SS(Y)_B \to SS(Z)_C$ *be soft mappings, then the composition of* $(f, e)$ *and* $(g, e')$ *is denoted by* $(g, e') \circ (f, e)$ *and defined by* $(g, e') \circ (f, e) = (g \circ f, e' \circ e)$.

## 2.2  Soft topological spaces

**Definition 2.2.1.** *[3] Let* $X$ *be an initial universe set,* $A$ *a set of parameters, and* $T \subseteq SS(X)_A$. *We say that the family* $T$ *defines a soft topology on* $X$ *if the following axioms are true:*

- a. $0_A, 1_A \in T$.
- b. *If* $(G, A), (H, A) \in T$, *then* $(G, A) \sqcap (H, A) \in T$.
- c. *If* $(G_i, A) \in T$ *for every* $i \in I$, *then* $\sqcup \{(G_i, A) : i \in I\} \in T$.

*The triplet* $(X, T, A)$ *is called a soft topological space. The members of* $T$ *are called soft open sets in* $X$. *Also, a soft set* $(F, A)$ *is called soft closed if the complement* $(F, A)^c$ *belongs to* $T$. *The family of soft closed sets is denoted by* $T^c$.

**Definition 2.2.2.** *[12] Let* $(X, T, A)$ *be a soft topological space and* $Y$ *be a nonempty subset of* $X$. *Then* $T_Y = \{(F, A)^Y : (F, A) \in T\}$, *is said to be the soft relative topology on* $Y$ *and* $(Y, T_Y, A)$ *is called a soft subspace of* $(X, T, A)$.

**Definition 2.2.3.** *[10] Let* $(X, T, A)$ *be a soft topological space over* $X$, $(F, A) \in SS(X)_A$ *and* $x \in X$. *Then* $(F, A)$ *is said to be a soft neighbourhood of* $x$ *if there exists a soft open set* $(G, A)$ *such that* $x \in (G, A) \sqsubseteq (F, A)$. *We denote the collection of all soft neighbourhoods of* $x$ *by* $\mathcal{U}_x$.

**Definition 2.2.4.** *[3] Let* $(X, T, A)$ *be a soft topological space. The soft closure* $\overline{(F, A)}$ *of* $(F, A) \in SS(X)_A$ *is the soft set* $\sqcap \{(G, A) \in SS(X)_A : (G, A) \in T^c, (F, A) \sqsubseteq (G, A)\}$.

**Proposition 2.2.5.** *[11] Let* $(X, T, A)$ *be a soft topological space,* $(F, A) \in SS(X)_A$ *and* $x \in X$. *Then* $x \in \overline{(F, A)}$ *implies that* $(F, A) \sqcap (G, A) \neq 0_A$, *for each* $(G, A) \in \mathcal{U}_x$.

### 2.2.1  Soft Separation Axioms

In this subsection, we introduce some separation axioms on **STop**.

**Definition 2.2.1.1.** *[12] Let $(X, T, A)$ be a soft topological space. Then,*

a. *$(X, T, A)$ is said to be soft $T_0$ if, whenever $x$ and $y$ are distinct elements of $X$, there is a soft open set $(G, A)$ containing one and not the other.*
b. *$(X, T, A)$ is said to be soft $T_1$ if, whenever $x$ and $y$ are distinct elements of $X$, there is a soft open set of each not containing the other.*
c. *$(X, T, A)$ is said to be soft $T_2$ or soft Hausdorff if, whenever $x$ and $y$ are distinct elements of $X$, there are soft open sets $(F, A)$ and $(G, A)$ such that $x \in (F, A)$, $y \in (G, A)$ and $(F, A) \sqcap (G, A) = 0_A$.*
d. *$(X, T, A)$ is said to be soft regular if, for each soft closed set $(F, A) \in SS(X)_A$ and each $x \in X$ such that $x \notin (F, A)$, there exist soft open sets $(G, A)$ and $(H, A)$ such that $x \in (G, A)$, $(F, A) \sqsubseteq (H, A)$ and $(G, A) \sqcap (H, A) = 0_A$.*
e. *$(X, T, A)$ is said to be soft $T_3$ if it is soft regular and soft $T_1$.*

**Definition 2.2.1.2.** *Let $x, y \in X$. We shall write $x \approx y$ to denote that $x \in (G, A)$, for each $(G, A) \in \mathcal{U}_y$ and $y \in (F, A)$, for each $(F, A) \in \mathcal{U}_x$.*

**Proposition 2.2.1.3.** *Let $(X, T, A)$ be a soft topological space. $(X, T, A)$ is soft $T_0$ iff for each $x, y \in X$, $x \approx y$ implies that $x = y$.*

*Proof:* If $(X, T, A)$ is soft $T_0$ and $x \neq y$, then there is a soft open set $(G, A)$ such that $x \in (G, A)$ and $y \notin (G, A)$ or there is an open set $(H, A)$ such that $y \in (H, A)$ and $\notin (H, A)$. In either case $x \not\approx y$. Conversely, if $(X, T, A)$ is not soft $T_0$. Then $x \neq y$ implies that, $y \in (F, A)$ and $x \in (G, A)$, for each $(F, A) \in \mathcal{U}_x$ and each $(G, A) \in \mathcal{U}_y$. This contradicts the hypothesis. ●

**Definition 2.2.1.4.** *A soft topological space $(X, T, A)$ is said to be soft $T_{0U}$ if for each $x, y \in X$, $x \approx y$ implies that, for each $a \in A$ and each $(F, A) \in T$, if $x \in F(a)$, then $y \in F(a)$ and for each $b \in A$ and each $(G, A) \in T$, if $y \in G(b)$, then $x \in G(b)$.*

**Proposition 2.2.1.5.** *Every soft $T_0$ space is soft $T_{0U}$.*

*Proof:* Suppose that $(X, T, A)$ is soft $T_0$, and choose $x, y \in X$ such that $x \approx y$. Let $a \in A$ and $(F, A) \in T$, such that $x \in F(a)$. Because $(X, T, A)$ is soft $T_0$, it follows that $x = y$. Thus $y \in F(a)$. Similarly, we can easily show that for each $b \in A$ and each $(G, A) \in T$ such that $y \in G(a)$, we have $x \in G(a)$.●

The converse of the preceding proposition doesn't always hold in general.

**Example 1.** Let $X = \{x, y, z\}$, $A = \{a_1, a_2\}$ and $T = \{0_A, 1_A, (F, A), (G, A)\}$ where, $F(a_1) = \{z\}$; $F(a_2) = \{z\}$; $G(a_1) = \{x, y\}$; $G(a_2) = \{x, y\}$. Clearly, $(X, T, A)$ is soft $T_{0U}$ but not soft $T_0$.

**Definition 2.2.1.6.** *Let $(X, T, A)$ be a soft topological space. Then,*

a. $(X,T,A)$ is said to be soft $T_{0^k}$ if, for each distinct elements $x, y \in X$, there is a soft open set $(F, A)$ such that $x \in (F, A)$ and $y \in (F, A)^c$, or there is a soft open set $(G, A)$ such that $y \in (G, A)$ and $x \in (G, A)^c$.
b. $(X,T,A)$ is said to be soft $T_{1^k}$ if, for each distinct elements $x, y \in X$, there are soft open sets $(F, A)$ and $(G, A)$ such that $x \in (F, A)$, $y \in (F, A)^c$, $y \in (G, A)$ and $x \in (G, A)^c$.
c. A soft topological space $(X, T, A)$ is said to be soft $k$ – regular if, for each soft closed set $(F, A) \in SS(X)_A$ and each $x \in X$ such that $x \in (F, A)^c$, there exist soft open sets $(G, A)$ and $(H, A)$ such that $x \in (G, A)$, $(F, A) \sqsubseteq (H, A)$ and $(G, A) \sqcap (H, A) = 0_A$.
d. $(X, T, A)$ is said to be soft $T_{3^k}$ if it is soft $k$ – regular and soft $T_{1^k}$.

**Proposition 2.2.1.7.** *Every soft $T_{0^k}$ space is soft $T_0$.*

*Proof:* The proof follows from the definition of soft $T_{0^k}$. ●

The converse of the preceding result doesn't hold in general.

**Example 2.** Let $X = \{x, y, z\}$, $A = \{a_1, a_2\}$ and $T = \{0_A, 1_A, (B, A), (C, A), (D, A)\}$, where, $B(a_1) = \{x\}$; $B(a_2) = \{x\}$; $C(a_1) = \{x, y\}$; $C(a_2) = \{y, z\}$; $D(a_1) = \{x\}$; $D(a_2) = \emptyset$. $(X, T, A)$ is clearly soft $T_0$ but not soft $T_{0^k}$.

Soft $T_{0^k}$ spaces need not be soft $T_1$, and soft $T_1$ spaces need not be soft $T_{0^k}$, as shown by the following examples.

**Example 3.** Let $X = \{x, y, z\}$, $A = \{a_1, a_2\}$ and $T = \{0_A, 1_A, (B, A), (C, A), (D, A), (E, A)\}$, where, $B(a_1) = \{x\}$; $B(a_2) = \{x\}$; $C(a_1) = \{x, y\}$; $C(a_2) = \{y\}$; $D(a_1) = \{x\}$; $D(a_2) = \emptyset$; $E(a_1) = \{x, y\}$; $E(a_2) = \{x, y\}$. Then, $(X, T, A)$ is soft $T_{0^k}$ but not soft $T_1$.

**Example 4.** Let $X = \{x, y, z\}$, $A = \{a_1, a_2\}$ and $T = \{0_A, 1_A, (B, A), (C, A), (D, A), (E, A), (F, A), (G, A), (H, A)\}$, where $B(a_1) = X$; $B(a_2) = \{x\}$; $C(a_1) = X$; $C(a_2) = \{y\}$; $D(a_1) = X$; $D(a_2) = \{z\}$; $E(a_1) = X$; $E(a_2) = \{x, y\}$; $F(a_1) = X$; $F(a_2) = \{x, z\}$; $G(a_1) = X$; $G(a_2) = \{y, z\}$; $H(a_1) = X$; $H(a_2) = \emptyset$. Then, $(X, T, A)$ is soft $T_1$ but not soft $T_{0^k}$.

**Proposition 2.2.1.8.** *Every soft $T_{1^k}$ space is soft $T_1$.*

*Proof:* Proof follows from the definition of soft $T_{1^k}$ space. ●

Example 4 shows that not every soft $T_1$ space is soft $T_{1^k}$.

**Proposition 2.2.1.9.** *A soft topological space $(X, T, A)$ is soft $T_{1^k}$ iff, for each $x \in X$, $(F_x, A)$ is soft closed.*

*Proof:* If $(X, T, A)$ is soft $T_{1^k}$ and $x \in X$, then each $y \neq x$ has a soft neighbourhood $(F, A)$ such that $(F_x, A) \sqcap (F, A) = 0_A$. Therefore, $(F_x, A)^c$ is a soft open set. Thus $(F_x, A)$ is soft closed. Conversely, if for each $x \in X$, $(F_x, A)$ is soft closed, then each $y \neq x$ belongs to $(F_x, A)^c$ and $x \in (F_y, A)^c$. ●

**Proposition 2.2.1.10.** *Every soft $T_2$ space is soft $T_{1^k}$.*

*Proof:* If $(X, T, A)$ is soft $T_2$ and $x, y \in X$ such that $x \neq y$. Then $x \in (F, A)$, $y \in (G, A)$ for some $(F, A), (G, A) \in T$ such that $(F, A) \sqcap (G, A) = 0_A$. Therefore $x \in (G, A)^c$ and $y \in (F, A)^c$. ●

The converse of the preceding result doesn't hold in general.

**Example 5.** Let $X = \mathbb{N}$, $A$ be a nonempty set and define soft open sets on $X$ as follows: $(F, A) \in T$ $iff$ $F(a) = \emptyset$ or $X \backslash F(a)$ $is$ $finite$, $for$ $every$ $a \in A$. Clearly, $(X, T, A)$ is soft $T_{1^k}$ but not soft $T_2$.

In Theorem 3.21 of [8], it was proved that every soft $T_3$ space is soft $T_{1^k}$. Using this result, we get that each soft $T_3$ space is soft $T_{3^k}$.

**Proposition 2.2.1.11.** *Every soft $T_3$ space is soft $T_{3^k}$.*

*Proof:* The proof relies on the fact that $x \in (F, A)^c$ implies $x \notin (F, A)$, for each $x \in X$ and each $(F, A) \in SS(X)_A$. ●

From the preceding result, we deduce that each soft regular space is soft $k$ −regular. However, the converse doesn't always hold, as shown by the following example.

**Example 6.** Let $X = \{x, y, z\}$, $A = \{a_1, a_2\}$ and $T = \{0_A, 1_A, (B, A), (C, A), (D, A), (E, A), (F, A), (G, A), (H, A), (I, A), (J, A), (K, A), (L, A), (M, A), (N, A), (O, A), (P, A), (Q, A)\}$, where, $B(a_1) = \{x\}$; $B(a_2) = \{x\}$; $C(a_1) = \{y\}$; $C(a_2) = \{y\}$; $D(a_1) = \{z\}$; $D(a_2) = \{z\}$; $E(a_1) = \{x\}$; $E(a_2) = \{y\}$; $F(a_1) = \{x, y\}$; $F(a_2) = \{x, y\}$; $G(a_1) = \{x, z\}$; $G(a_2) = \{x, z\}$; $H(a_1) = \{x\}$; $H(a_2) = \{x, y\}$; $I(a_1) = \{y, z\}$; $I(a_2) = \{y, z\}$; $J(a_1) = \{x, y\}$; $J(a_2) = \{y\}$; $K(a_1) = \{x, z\}$; $K(a_2) = \{y, z\}$; $L(a_1) = \{x, z\}$; $L(a_2) = X$; $M(a_1) = \{z\}$; $M(a_2) = \{y, z\}$; $N(a_1) = \{x\}$; $N(a_2) = \emptyset$; $O(a_1) = \emptyset$; $O(a_2) = \{y\}$; $P(a_1) = X$; $P(a_2) = \{y, z\}$; $Q(a_1) = \{x, z\}$; $Q(a_2) = \{z\}$. Clearly, $(X, T, A)$ is soft $T_{3^k}$ but not soft $T_3$.

**Proposition 2.2.1.12.** *Every soft $T_{3^k}$ space is soft $T_2$.*

*Proof:* Suppose that $(X, T, A)$ is a soft $T_3^k$ space. Let $x, y \in X$ such that $x \neq y$. Then there are soft open sets $(F, A)$ and $(G, A)$ in $X$ such that $x \in (F, A)$, $y \in (F, A)^c$, $y \in (G, A)$ and $x \in (G, A)^c$. But $(X, T, A)$ is soft $k$ −regular and $(F, A)^c$ is soft closed, so there are soft open sets $(H, A)$ and $(J, A)$ such that $x \in (H, A)$, $(F, A)^c \sqsubseteq (J, A)$ and $(H, A) \sqcap (J, A) = 0_A$. Thus $(X, T, A)$ is soft $T_2$. ●

**Proposition 2.2.1.13.** *Let $(X, T, A)$ be a soft $T_\alpha$ space, where $\alpha = 0^k, 1^k, 3^k$, and $Y \subseteq X$ be a nonempty set. Then $(Y, T_Y, A)$ is soft $T_\alpha$.*

*Proof:* We only prove for $\alpha = 0^k$, as well as soft $k$-regularity. If $(X, T, A)$ is soft $T_{0^k}$ and $x, y \in Y$ such that $x \neq y$, then there is a soft open set $(F, A)$ in $X$ such that $x \in (F, A)$ and $y \in (F, A)^c$ or there is a soft open set $(G, A)$ in $X$ such that $y \in (G, A)$ and $x \in (G, A)^c$. We only show the first case: $x \in (F, A)$ and $y \in (F, A)^c$. We have that $x \in 1_A^Y \sqcap (F, A) = (F, A)^Y$ and $y \in 1_A^Y \setminus (F, A)^Y$, respectively. Thus $(Y, T_Y, A)$ is soft $T_{0^k}$.

For soft $k$-regularity: Let $(F, A)^Y$ be a soft closed set of $Y$ and let $x \in 1_A^Y \setminus (F, A)^Y$. Since $(F, A)^Y = 1_A^Y \sqcap (G, A)$, for some $(G, A) \in T^c$, we have that $x \in 1_A \setminus (G, A)$. Therefore, there are soft open sets $(H, A)$ and $(K, A)$ in $X$ such that $x \in (H, A)$, $(G, A) \sqsubseteq K, A)$ and $(H, A) \sqcap (K, A) = 0_A$. Let $(J, A) = 1_A^Y \sqcap (H, A)$ and $(L, A) = 1_A^Y \sqcap (J, A)$. Then $(J, A), (L, A) \in T_Y$, $x \in (J, A)$ and $(J, A) \sqcap (L, A) = 0_A$. Thus $(Y, T_Y, A)$ is soft $k$-regular. ●

**Proposition 2.2.1.14.** *The property of being soft $T_\alpha$, where $\alpha = 0^k, 1^k, 3^k$, is preserved under a soft homeomorphism mapping.*

*Proof:* Suppose that $g = (f, e): (X, T, A) \to (X', T', B)$ is a soft mapping such that $g$ is soft onto, soft injective, both $g$ and $g^{-1}$ are soft continuous. We prove for $\alpha = 0^k$ and soft $k$-regularity. Assume that $(X, T, A)$ is a soft $T_{0^k}$ space. Let $x, y \in X'$ such that $x \neq y$. Then there exist $x', y' \in X$ such that $f(x') = x$, $f(y') = y$ and $x' \neq y'$. But $(X, T, A)$ is soft $T_0^k$, so there is a soft open set $(F, A)$ such that $x' \in (F, A)$ and $y' \in (F, A)^c$ or there is a soft open set $(G, A)$ such that $y' \in (G, A)$ and $x' \in (G, A)^c$. We prove for the first case: $x' \in (F, A)$ and $y' \in (F, A)^c$. We have that $x \in g((F, A))$ and $y \in g((F, A)^c) = (g((F, A)))^c$, where, by soft continuity of $g^{-1}$, $g((F, A))$ is soft open. Thus $(X', T', B)$ is soft $T_{0^k}$.

Soft $k$-regularity: Suppose that $(X, T, A)$ is soft $k$-regular. Let $x \in X'$ and $(F, B) \in (T')^c$ such that $x \in (F, B)^c$. Then there exists $x' \in X$ such that $f(x') = x$ and $x' \in g^{-1}((F, B)^c)$, where, by soft continuity of $g$, $g^{-1}((F, B)^c) \in T^c$. But $(X, T, A)$ is soft $k$-regular, so there are soft sets $(G, A), (H, A) \in T$ such that $x' \in (G, A)$, $g^{-1}((F, B)^c) \sqsubseteq (H, A)$ and $(G, A) \sqcap (H, A) = 0_A$. It follows that $x \in g((G, A))$, $(F, B)^c \sqsubseteq g((H, A))$ and $g((G, A)) \sqcap g((H, A)) = 0_B$, where $g((G, A)), g((H, A)) \in T'$. Thus $(X', T', B)$ is soft $k$-regular.

We consider soft topological spaces whose points are separated by soft clopen (simultaneously soft closed and soft open) sets.

**Definition 2.2.1.15.** *A soft topological space $(X, T, A)$ is soft totally separated if, whenever $x$ and $y$ are distinct elements of $X$, there is a soft clopen set $(G, A)$ containing one and not the other.*

**Definition 2.2.1.16.** *A soft topological space $(X, T, A)$ is said to be soft $k$-totally separated if for each distinct elements $x, y \in X$, there is a soft clopen set $(F, A)$ such that $x \in (F, A)$ and $y \in (F, A)^c$, or there is a soft clopen set $(G, A)$ such that $y \in (F, A)$ and $x \in (G, A)^c$.*

**Proposition 2.2.1.17.** *Every soft $k$-totally separated space is soft totally separated.*

*Proof:* Proof easily follows from the definition of a soft $k$-totally separated space. ●

The converse of the preceding proposition doesn't hold in general.

**Example 7.** Let $X = \{x, y\}$, $A = \{a_1, a_2\}$ and $T = \{0_A, 1_A, (B, A), (C, A)\}$, where $B(a_1) = \{x\}$; $B(a_2) = X$; $C(a_1) = \{y\}$; $C(a_2) = \emptyset$. Clearly, $(X, T, A)$ is soft totally separated but not soft $k$ −totally separated.

The preceding example also shows that there are some soft totally separated spaces which are neither soft $T_0^k$ nor soft $T_1$. However, every soft $k$ −totally separated space is soft $T_2$.

**Proposition 2.2.1.18.** *Every soft $k$ −totally separated space is soft $T_2$.*

*Proof:* Suppose that $(X, T, A)$ is a soft $k$ −totally separated space. Let $x, y \in X$ such that $x \neq y$. Then there is a soft clopen set $(F, A)$ such that $x \in (F, A)$ and $y \in (F, A)^c$ or there is a soft clopen set $(G, A)$ such that $y \in (G, A)$ and $x \in (G, A)^c$. Consider $x \in (F, A)$ and $y \in (F, A)^c$. We have that $(F, A) \sqcap (F, A)^c = 0_A$ and $(F, A)^c \in T$. Thus $(X, T, A)$ is soft $T_2$. Similarly, we can show that $y \in (G, A)$ and $x \in (G, A)^c$ implies that $(X, T, A)$ is soft $T_2$. ●

Example 6 shows that not every soft $k$ −totally separated space is soft $T_3$.

**Proposition 2.2.1.19.** *Every soft $k$ −totally separated space is soft $T_{3^k}$.*

*Proof:* It suffices to show that a soft $k$ −totally separated space $(X, T, A)$ is soft $k$ −regular. If $x \in (F, A)^c$ for some $(F, A)^c \in T^c$, then for each $y \in (F, A)$, there is a soft clopen set $(G, A)_y$ such that $x \in (G, A)_y$ and $y \in (G, A)_y^c$. We have that $(F, A) \sqsubseteq \sqcup\{(G, A)_y^c : y \in (F, A)\}$. Therefore $\sqcup\{(G, A)_y^c : y \in (F, A)\}$ and $(G, A)_y$ are the required soft open sets. ●

The following result follows from the definition of soft totally separated spaces. Its proof shall be omitted.

**Proposition 2.2.1.20.** *Every soft totally separated space is soft $T_0$.*

**Proposition 2.2.1.21.** *Let $(X, T, A)$ be a soft topological space and $Y \subseteq X$ be a nonempty set. Then $(Y, T_Y, A)$ is soft totally separated (respectively, soft $k$ −totally separated) whenever $(X, T, A)$ is soft totally separated (respectively, soft $k$ −totally separated).*

*Proof:* The proof is similar to the first part of the proof of Proposition 2.2.1.13. ●

**Proposition 2.2.1.22.** *The property of being soft totally separated ($k$ −totally separated) is preserved under a soft homeomorphism mapping.*

*Proof:* The proof is similar to the first part of the proof of Proposition 2.2.1.14. ●

# 3. Soft $T_0$ Reflection

Denote by $\mathbf{STop_0}$ the full subcategory of $\mathbf{STop}$ whose objects are soft $T_0$ spaces and morphisms are soft continuous mappings. For the rest of this text, we assume that all soft topological spaces belong to $\mathbf{STop_{0U}}$. The aim of this section is to construct the soft $T_0$ reflection of a soft $T_{0U}$ topological space and find some of its properties that shall be useful in the next section.

Consider the equivalence relation on $(X, T, A)$ given by: $x \sim y$ $iff$ $x \approx y$. Denote, by $X_0$ and $r_X$, the collection of all $\sim$-equivalence classes on $(X, T, A)$, and the canonical surjection from $X$ to $X_0$, respectively. Consider the identity map $i: A \to A$. Clearly, the map $g_{0X} = (r_X, i): SS(X)_A \to SS(X_0)_A$ is a soft onto mapping. Define a soft topology on $X_0$ by $T_0 = \{(F, A) \in SS(X_0, A): g_{0X}^{-1}((F, A)) \in T\}$. It can be easily seen that $(X_0, T_0, A)$ is a soft topological space and $g_{0X}: SS(X)_A \to SS(X_0)_A$ is a soft continuous and soft onto map.

**Lemma 3.1.** *Suppose that $(X, T, A)$ is a soft topological space. Then $(F, A) = g_{0X}^{-1}(g_{0X}((F, A)))$, for each $(F, A) \in T$.*

*Proof:* Let $a \in A$ and $(H, A) = g_{0X}^{-1}(g_{0X}((F, A)))$. If $x \in H(a)$, then $[x] \in g((H, A))$. Let $(G, A) = g((H, A))$. We have that $[x] \in G(a)$. But $r_X$ is onto, so $[x] = [y]$ for some $y \in F(a)$. Therefore $x \approx y$. Because $(X, T, A)$ is soft $T_{0U}$, it follows that $x \in F(a)$. Thus $g_{0X}^{-1}(g_{0X}((F, A))) \sqsubseteq (F, A)$. ●

**Lemma 3.2.** *The soft map $g_{0X}$ is soft open.*

*Proof:* Follows from the combination of the preceding result and the fact that $g_{0X}$ is soft onto. ●

**Lemma 3.3.** *$(X_0, T_0, A)$ is soft $T_0$.*

*Proof:* Follows from the definition of $x \approx y$ and the preceding lemma. ●

**Theorem 3.4.** *The pair $(g_{0X}, (X_0, T_0, A))$ is the soft $T_0$ reflection of $(X, T, A)$.*

*Proof:* Let $f: X \to Y$ and $e: A \to B$ be functions such that $(f, e): (X, T, A) \to (Y, T', B)$ is a soft continuous mapping, where $(Y, T', B)$ is soft $T_0$. We have that $x \sim y \Rightarrow f(x) \approx f(y)$. Because $Y$ is soft $T_0$, we have that $f(x) = f(y)$. Define $f': X_0 \to Y$ by $f'(r_X(x)) = f(x)$. Then $(f', e)$ is a unique soft continuous and soft open map satisfying $(f', e) \circ g_{0X} = (f, e)$. ●

**Corollary 3.5.** *$\mathbf{STop_0}$ is a reflective subcategory of $\mathbf{STop_{0U}}$.*

**Definition 3.6.** *Let $(X, T, A)$ and $(Y, T', B)$ be soft topological spaces. A soft continuous map $(f, e): SS(X)_A \to SS(Y)_B$ is said to be soft initial if for each $(F, A) \in T, (F, A) = (f, e)^{-1}((G, B))$ for some $(G, B) \in T'$.*

**Proposition 3.7.** Let $(X,T,A)$ and $(Y,T',B)$ be soft topological spaces and $(f,e):(X,T,A) \to (Y,T',B)$ be a soft continuous mapping. Then $(f,e)$ is soft initial if and only if, for each soft closed set $(F,A) \in SS(X)_A$, $(F,A) = (f,e)^{-1}(\overline{(f,e)((F,A))})$.

*Proof:* It suffices to show that, for each $(F,A) \in T^c$, there is $(G,B) \in T'$ such that $(f,e)^{-1}((G,B)^c) = (f,e)^{-1}(\overline{(f,e)((F,A))})$. ●

**Definition 3.8.** A soft set $(F,A) \in SS(X)_A$ is said to be soft locally closed if it can be written as a soft intersection of a soft closed set and a soft open set.

**Definition 3.9.** Let $(X,T,A)$ and $(Y,T',B)$ be soft topological spaces. A soft continuous map $(f,e):SS(X)_A \to SS(Y)_B$ is said to be a soft quasihomeomorphism if it is soft initial and $(f,e)(1_A)$ meets each soft locally closed set in $(Y,T',B)$.

For soft quasihomeomorphisms, we get the following result which has a straight forward proof.

**Proposition 3.10.** Let $(X,T,A), (Y,T',B)$ and $(Z,T'',C)$ be soft topological spaces. Let $(f,e):(X,T,A) \to (Y,T',B)$ and $(h,e'):(Y,T',B) \to (Z,T'',C)$ be soft continuous mappings. Then, if two among $(f,e), (h,e')$ and $(h,e') \circ (f,e)$ are soft quasihomeomorphisms, then so is the third one.

**Proposition 3.11.** Let $(X,T,A)$ and $(Y,T',B)$ be soft topological spaces and $(f,e):SS(X)_A \to SS(Y)_B$ be a soft quasihomeomorphism. Then the following statements hold:

a. If $(X,T,A)$ is soft $T_0$ and $e$ is injective, then $(f,e)$ is injective.
b. If $(Y,T',B)$ is soft $T_{1^k}$ and $e$ is surjective, then $(f,e)$ is onto.
c. If $(X,T,A)$ is soft $T_0$ and $(Y,T',B)$ is soft $T_{1^k}$ and $e$ is bijective, then $(f,e)$ is a soft homeomorphism.

*Proof:* We only prove (a). If $x,y \in X$ such that $x \neq y$, then there is a soft open set $(F,A)$ such that $x \in (F,A)$ and $y \notin (F,A)$. Because $(f,e)$ is a soft quasihomeomorphism, then $(F,A) = (f,e)^{-1}((G,B))$, for some $(G,B) \in T'$. Therefore $x \in (G,B)$ and $y \notin (G,B)$. Thus $f(x) \neq f(y)$. ●

**Proposition 3.12.** $g_{0X}$ is a soft quasihomeomorphism.

*Proof:* Follows since $(F,A) = g_{0X}^{-1}(g_{0X}((F,A)))$, for all $(F,A) \in T$ and the fact that $g_{0X}$ is soft onto. ●

**Proposition 3.13.** $g_{0X}$ is soft closed.

*Proof:* Follows from Lemma 3.2 and the fact that $g_{0X}((F,A)^c) = \big(g_{0X}((F,A))\big)^c$, for all $(F,A) \in SS(X)_A$. ●

In what follows, we denote by $(T_0(f), e)$, the soft map between $(X_0, T_0, A)$ and $(Y_0, T'_0, B)$, where $T_0(f)$ is defined as $r_X(x) = r_Y(f(x))$. It can easily be shown that this soft map is soft continuous and $g_{0Y} \circ (f, e) = (T_0(f), e) \circ g_{0X}$.

**Proposition 3.14.** *Let $(X, T, A)$ and $(Y, T', B)$ be soft topological spaces and $(f, e): SS(X)_A \rightarrow SS(Y)_B$ be a soft continuous map. Then $(f, e)$ is a soft quasihomeomorphism iff $(T_0(f), e)$ is a soft quasihomeomorphism.*

*Proof:* Follows from Proposition 3.10. ●

## 4. Some New Separation Axioms

In this section, we introduce and characterize soft $T_{(0,\alpha)}$ spaces.

**Definition 4.1.** *Let $\alpha \in \{0^k, 1, 1^k, 2, 3^k, 3, TS, (TS)^k\}$. A soft topological space $(X, T, A)$ is said to be a soft $T_{(0,\alpha)}$ space if the soft $T_0$ reflection of $(X, T, A)$ is a soft $T_{0^k}$ (respectively, soft $T_1$, soft $T_{1^k}$, soft $T_2$, soft $T_{3^k}$, soft $T_3$, soft totally separated, soft $k$ − totally separated)-space.*

Some of these separation axioms were first defined in Definition 1.1 of [2] on **Top**. They were then studied in different categories such as **ORDTOP** (category of ordered topological spaces whose morphisms are continuous increasing maps) in [6], **GenTop** in [9] and **PreTop** in [7].

**Proposition 4.2.** *Let $(X, T, A)$ be a soft topological space. Then $(X, T, A)$ is a soft $T_{(0,0^k)}$ space iff for each $x, y \in X$ such that $x \not\approx y$, there exist a soft open set $(G, A)$ such that $x \in (G, A)$ and $y \in (G, A)^c$ or there is a soft open set $(F, A)$ such that $y \in (F, A)$ and $x \in (F, A)^c$.*

*Proof:* If $x, y \in X$ such that $x \not\approx y$, then $[x] \neq [y]$. Since $(X_0, T_0, A)$ is soft $T_{0^k}$, there is a soft open set $(\mathcal{H}, A)$ in $X_0$ such that $[x] \in (\mathcal{H}, A)$ and $[y] \in (\mathcal{H}, A)^c$ or there is a soft open set $(\mathcal{K}, A)$ in $X_0$ such that $[y] \in (\mathcal{K}, A)$ and $[x] \in (\mathcal{K}, A)^c$. We only prove the first case: $[x] \in (\mathcal{H}, A)$ and $[y] \in (\mathcal{H}, A)^c$. We have that $y \in \left(g_{0X}^{-1}((\mathcal{H}, A))\right)^c$, where $g_{0X}^{-1}((\mathcal{H}, A)) \in T$ and $x \in g_{0X}^{-1}((\mathcal{H}, A))$.
Conversely, choose $x, y \in X$ such that $[x] \neq [y]$. Then $x \not\approx y$. By hypothesis, we have $(G, A) \in T$ such that $x \in (G, A)$ and $y \in (G, A)^c$ or $(H, A) \in T$ such that $y \in (H, A)$ and $x \in (H, A)^c$. We only consider the first case: $x \in (G, A)$ and $y \in (G, A)^c$. Since $g_{0X}$ is a soft quasihomeomorphism, there is a soft open set $(\mathcal{H}, A)$ in $(X_0, T_0, A)$ such that $(G, A) = g_{0X}^{-1}((\mathcal{H}, A))$. It is clear that $[y] \in (\mathcal{H}, A)^c$ and $[x] \in (\mathcal{H}, A)$. Thus $(X_0, T_0, A)$ is soft $T_{0^k}$. ●

**Proposition 4.3.** *Let $(X, T, A)$ be a soft topological space. Then $(X, T, A)$ is a soft $T_{(0,1)}$ space iff, for each $x, y \in X$ such that $x \not\approx y$, there are soft open sets $(G, A)$ and $(H, A)$ in $X$ such that $x \in (G, A), y \in (H, A), x \notin (H, A)$ and $y \notin (G, A)$.*

*Proof:* The proof is similar to the proof of the preceding result. ●

**Proposition 4.4.** Let $(X, T, A)$ be a soft topological space. Then the following statements are equivalent.

a. $(X, T, A)$ is a soft $T_{(0,1^k)}$ space.
b. For each $x, y \in X$ such that $x \not\approx y$, there are soft open sets $(G, A)$ and $(H, A)$ such that $x \in (G, A)$, $y \in (H, A)$, $x \in (H, A)^c$ and $y \in (G, A)^c$.
c. For each $x \in X$, the soft set $(F_{\dot{x}}, A)$ defined as $F_{\dot{x}}(a) = [x]$, for all $a \in A$, is soft closed in $X$.

*Proof:* To show that $a \Leftrightarrow c$, it suffices to show that $(F_{\dot{x}}, A) = g_{0X}^{-1}((F_{[x]}, A))$ and $g_{0X}((F_{\dot{x}}, A)) = (F_{[x]}, A)$ for all $x \in X$. ●

**Proposition 4.5.** Let $(X, T, A)$ be a soft topological space. Then $(X, T, A)$ is a soft $T_{(0,2)}$ space iff for each $x, y \in X$ such that $x \not\approx y$, there are soft open sets $(G, A)$ and $(H, A)$ such that $x \in (G, A)$, $y \in (H, A)$ and $(G, A) \sqcap (H, A) = 0_A$.

**Proposition 4.6.** Let $(X, T, A)$ be a soft topological space. Then the following statements are equivalent:

a. $(X, T, A)$ is a soft $T_{(0,3^k)}$ space.
b. The following statements hold:
   i. $(X, T, A)$ is a soft $T_{(0,1^k)}$ space.
   ii. $(X, T, A)$ is a soft $k$-regular space.

*Proof:* The proof follows since $g_{0X}$ is soft onto, soft closed and soft initial. ●

**Proposition 4.7.** Let $(X, T, A)$ be a soft topological space. Then the following statements are equivalent:

a. $(X, T, A)$ is a soft $T_{(0,3)}$ space.
b. The following statements hold:
   i. $(X, T, A)$ is a soft $T_{(0,1)}$ space.
   ii. $(X, T, A)$ is a soft regular space.

*Proof:* The proof follows since $g_{0X}$ is soft onto, soft closed and soft initial. ●

**Proposition 4.8.** Let $(X, T, A)$ be a soft topological space. Then $(X, T, A)$ is a soft $T_{(0,TS)}$ space iff, for each $x, y \in X$ such that $x \not\approx y$, there exist a soft clopen set $(G, A)$ such that $x \in (G, A)$ and $y \notin (G, A)$ or there is a soft clopen set $(F, A)$ such that $y \in (F, A)$ and $x \notin (F, A)$.

**Proposition 4.9.** Let $(X, T, A)$ be a soft topological space. Then $(X, T, A)$ is a soft $T_{(0,(TS)^k)}$ space iff for each $x, y \in X$ such that $x \not\approx y$, there exist a soft clopen set $(G, A)$ such that $x \in (G, A)$ and $y \in (G, A)^c$ or there is a soft clopen set $(F, A)$ such that $y \in (F, A)$ and $x \in (F, A)^c$.

**Proposition 4.10.** *Let $(X, T, A)$ and $(Y, T', B)$ be soft topological spaces and suppose that $(f, e): SS(X)_A \to SS(Y)_B$ is a soft onto quasihomeomorphism. Then $(X, T, A)$ is a soft $T_{(0,\alpha)}$ space, where $\alpha = 0^k, 1, TS$, iff $(Y, T', B)$ is a soft $T_{(0,\alpha)}$ space.*

*Proof:* This follows since $(T_0(f), e)$ is a soft homeomorphism. ●

**Proposition 4.11.** *Let $(X, T, A)$ and $(Y, T', B)$ be soft topological spaces and suppose that $(f, e): SS(X)_A \to SS(Y)_B$ is a soft quasihomeomorphism. Then $(X, T, A)$ is a soft $T_{(0,\alpha)}$ space, where $\alpha = 1^k, 2, 3^k, 3, (TS)^k$ iff $(Y, T', B)$ is a soft $T_{(0,\alpha)}$ space*

*Proof:* This follows since $(T_0(f), e)$ is a soft homeomorphism. ●

**Proposition 4.12.** *Every soft subspace of a soft $T_{(0,\alpha)}$ space is soft $T_{(0,\alpha)}$.*

*Proof:* Follows since every soft subspace of a soft $T_\alpha$ space is soft soft $T_\alpha$.

**Example 8.** The soft topological space given in Example 1 is an example of a soft $T_{(0,\alpha)}$ space, where $X_0 = \{[x], [z]\}, T_0 = \{(H, A), (J, A), 1_A^{X_0}, 0_A\}$ and $H(a_1) = \{[x]\}$; $H(a_2) = \{[x]\}$; $J(a_1) = \{[z]\}$; $J(a_2) = \{[z]\}$.

**Example 9.** Let $X = \{x, y, z\}, A = \{a_1, a_2\}$, and $T = \{0_A, 1_A, (B, A), (C, A)\}$, where $B(a_1) = \{x\}$; $C(a_1) = \{x\}$; $B(a_2) = \emptyset$; $C(a_2) = \{x\}$. Then, $(X, T, A)$ is a soft $T_{0U}$ space that is soft $T_{(0,0^k)}$ but neither soft $T_{(0,1)}$ nor soft $T_{(0,TS)}$. This means that $(X, T, A)$ is not soft $T_{(0,\alpha)}$, for $\alpha \in \{1, 1^k, 2, 3^k, 3, (TS)^k\}$.